\documentclass[11pt,reqno]{amsart}
\usepackage{amsmath,bm}
\usepackage{amsfonts}
\usepackage{amssymb}
\usepackage{amsthm}
\usepackage{amscd}
\usepackage{a4wide}
\usepackage{tikz}\usepackage{amsmath}
\usepackage{amssymb}
\usepackage{amsthm}
\usepackage{graphics}
\usepackage{eucal}
\usepackage{mathrsfs}
\usepackage[utf8]{inputenc} % Any characters can be typed directly from the keyboard, eg éçñ
\usepackage{textcomp} % provide lots of new symbols
\usepackage{graphicx}  % Add graphics capabilities
\usepackage{flafter}  % Don't place floats before their definition
\usepackage{bm}  % Define \bm{} to use bold math fonts

\usepackage[pdftex,bookmarks,colorlinks,breaklinks]{hyperref}  % PDF hyperlinks, with coloured links
\usepackage{color}
\usepackage{memhfixc}  % remove conflict between the memoir class & hyperref
\usepackage{pdfsync}  % enable tex source and pdf output syncronicity
\usepackage{amsmath,amsthm,indentfirst,marvosym}
\usepackage{amssymb}
\usepackage{amsfonts}
\theoremstyle{plain}
\newtheorem*{maintheorem*}{Main Theorem}
\newtheorem*{thm*}{Theorem}
\newtheorem*{thma*}{Theorem A}
\newtheorem*{thmaa*}{Theorem A'}
\newtheorem*{thmb*}{Theorem B}
\newtheorem*{thmo*}{Theorem 1.1}
\newtheorem*{thmc*}{Theorem C}
\newtheorem*{thmd*}{Theorem D}
\newtheorem*{thmf*}{Theorem 4.1}
\newtheorem*{remark*}{Remark}
\newtheorem*{conjecture*}{Conjecture}
\newtheorem*{prop*}{Proposition}
\newtheorem*{lem*}{Basic Lemma}
\newtheorem{thm}{Theorem}[section]
\newtheorem{cor}[thm]{Corollary}
\newtheorem{lem}[thm]{Lemma}
\newtheorem{prop}[thm]{Proposition}

\theoremstyle{definition}

\newtheorem*{proofc*}{Proof of Theorem C}

\newtheorem{conjecture}[thm]{Conjecture}

\newtheorem{definition}[thm]{Definition}

\def\Bc{\mathcal{B}}
\def\bbz{\mathbb{Z}}
\def\bbq{\mathbb{Q}}

\def\bbr{\mathbb{R}}

\def\bbc{\mathbb{C}}

\def\gcal{\mathcal{G}}

\def\SL{\rm{SL}}

\def\V{\mathcal{V}}

\def\hh{\hspace{0.5mm}}

\begin{document}

%\thanks {\textit{Keywords: Diophantine approximation, unipotent dynamics, quadratic forms}}
\author{Youssef Lazar}
\address{{\tiny Y. Lazar,  School of Mathematics, University of East Anglia, Norwich NR4 7TJ UK}}
\email{ ylazar77@gmail.com}
\date{}
\title[$S$-arithmetic Oppenheim conjecture for pairs $(Q,L)$]{ Values of pairs involving one quadratic form and one linear form at $S$-integral points}

%%%%%%%%%%%%%%%%%%%%% Abstract %%%%%%%%%%%%%%%%%%%%%%%%%%%%%%%%%%%%%%%%%%%%%%%%%%%%%%%%%%%%%%%%%%%%%%%%%%%
\maketitle

\begin{abstract} 
  We prove the existence of $S$-integral solutions of  simultaneous diophantine inequalities for pairs $(Q,L)$ involving one quadratic form and one linear form satisfying some arithmetico-geometric conditions. This result generalises previous results of Gorodnik and Borel-Prasad. The proof uses Ratner's theorem for unipotent actions on homogeneous spaces combined with an argument of strong approximation.  \end{abstract}

%%%%%%%%%%%%%%%%%%%%%%%%%%%%%%%%%%%%%%%%%%%%%%%%%%%%%%%%%%%%%%%%%%%%%%   Introduction   %%%%%%%%%%%%%%%%%%%%%%%%%%%%%%%%%%%%%%%%%%%%%%%%%%%%%%%%%%%%%%%%%%%%%%%%%%%
\section{Introduction}
\label{intro}

 The theory of unipotent flows on homogeneous spaces is a powerful tool used to solve many difficult problems in number theory and more particularly in diophantine approximation. One of the great achievement of those so-called dynamical methods is the proof made by G.A. Margulis of the Oppenheim conjecture: \textit{Let  $Q$ be a nondegenerate indefinite real quadratic form in $n \geqslant 3$ variables which is not proportional to a form with rational coefficients then $Q(\bbz^{n})$ is dense in $\bbr$}. A similar Oppenheim type problem concerns the existence of integral solutions of simultaneous diophantine inequalities involving one quadratic form and one linear form. More precisely given a pair $(Q,L)$ and $(a,b) \in \bbr^{2}$ the problem is to find sufficient conditions which guarantees the existence of a nonzero integral vector in $ x \in \bbz^{n}$ such that 
 \begin{center}
$(\mathrm{A})$  For any $ \varepsilon>0$ one has simultaneously  $  |Q(x) - a|<  \varepsilon $  and $     |L(x) - b|<  \varepsilon    $\end{center}  
This condition is equivalent to ask the density of the set  $\{ (Q(x), L(x)) : x \in \bbz^{n}\} $ in $\bbr^{2}$. The first result in that direction is due to S.G Dani and G.A. Margulis \cite{DM}  and concerns the dimension 3 for a pair $(Q,L)$ consisting of one nondegenerate indefinite quadratic form and a nonzero linear form in dimension $3$ such that the cone $\{Q = 0\}$ intersects tangentially the plane $\{L = 0\}$  and no linear combinaison of $Q$ and $L^{2}$ is rational. Under those conditions they proved using the original method used to prove the Oppenheim conjecture that the set $\{ (Q(x), L(x)) : x \in \bbz^{3}\} $ is dense in $\bbr^{2}$. In higher dimension,  the density for pairs holds if one replaces the previous transversality condition by the assumption that  $Q_{|L = 0}$ is indefinite, this result is due to A.Gorodnik \cite{G1}: 

\begin{thm}[Gorodnik]~\label{thm} Let $F = (Q,L)$ be a pair consisting of a quadratic form $Q$ and $L$ a nonzero linear form in dimension $n \geq 4$ satisfying the the following conditions
\begin{enumerate}
  \item $Q$ is nondegenerate.
  \item $Q_{|L = 0}$ is indefinite.
  \item  No linear combination of $Q$ and $L^{2}$ is rational.
  \end{enumerate}
Then the set $F(\mathcal{P}( \bbz^{n})) $ is dense in $\bbr^{2}$ where $\mathcal{P}( \bbz^{n})$ is the set of primitive integer vectors.
\end{thm}
\noindent  The conclusion of the theorem implies immediately that the set  $F( \bbz^{n})$ is dense in $\bbr^{2}$.
The proof of this theorem reduced to the case of the dimension $4$. The condition $(1)$ is a sufficient condition to  ensure that we have $ F(\bbr^{n})  =  \bbr^{2}$ and this is a conjecture that this condition can be weakened in order to make it necessary (see \S \ref{sec end}). The most important obstruction to prove density for pairs is that the identity component of the stabilizer of a pair $(Q,L)$ is no longer maximal among the connected Lie subgroups of $\SL(4,\bbr)$ in contrast with the case of the isotropy groups $\rm{SO}$$(3,1)^\circ$ or $\rm{SO}$$(2,2)^\circ$.
 
\noindent  The stabilizer of the pair $(Q,L)$ is defined by the following subgroup of $G$,
   \begin{center}
$\rm{Stab}$$(Q,L) = \big\{ h \in SO(Q)   |   L(hx)=L(x)   \}  $. 
\end{center} 

\noindent It is not difficult to see that there exists $g \in G$ such that  $ (Q,L) = (Q_{0}^{g}, L_{0}^{g} )$ for some canonical pairs $ (Q_{0},L_{0})$ given explicitly  (see \cite{G1}, Proposition 2). Clearly one has $\rm{Stab}$$(Q,L) = g$$\rm{Stab}$$(Q_{0},L_{0})g^{-1}$  and we are reduced to study the stabilizer of canonical pairs $ (Q_{0},L_{0})$. The pairs such that $Q_{| L=0}$ is nondegenerate (resp. degenerate) are said to be of type (I) (resp. II).
    The proof of Theorem \ref{thm} is divided in two parts following each type and consists to apply Ratner's orbit closure theorem, and to study the action of  $\rm{Stab}$$(Q_{0},L_{0})$ on the dual space of $\bbc^{4}$. A remarkable fact is that the density is proved without showing the density of the orbit closure of the stabilizer in the homogeneous space $G/\Gamma$.  Indeed the intermediate subgroups which possess non-trivial irreducible components have closed orbits in $G/\Gamma$, in particular they are not maximal. However, one is able to classify all the complex semisimple Lie algebras in $\mathfrak{sl}$$(4,\bbc)$, and Gorodnik used this classification to check density case by case using the constrain on rationality given by the condition (3). The situation for pairs of type (II) is more complicated compared with the pairs of type (I) since the dual action of the stabilizer has three irreducible components for the pairs of type (II), instead of two for the pairs of type (I). \\
 We are going to show an  $S$-arithmetic generalisation of this result for pairs of type (I). Our proof is influenced by the work of Borel-Prasad on generalised the Oppenheim conjecture for quadratic forms (\cite{BP}) and of course also by Gorodnik's  proof of  theorem~\ref{thm}.
 \medskip
  
   \section{ Main result } 

 \medskip
 
 \subsection{$S$-arithmetic setting} Let us recall what we mean by $S$-arithmetic setting by fixing some notations.
Let $k$ be a number field, that is a finite extension of $\bbq$ and let $\mathcal{O}$ be the ring of integers of $k$.
 For every normalised absolute value $\vert. \vert_{s}$ on $k$, let $k_s$ be the completion of $k$ at $s$. We identify $s$ with the specific absolute value $| .|_s $ on $k_s$ defined by the formula $\mu(a\Omega) = |a|_s \mu(\Omega)$, where $\mu$ is any Haar measure on the additive group $k_s$, $a \in k_s$  and $\Omega $ is a measurable subset of $k_s$ of finite measure. We denote by $\Sigma_k$ the set of places of $k$.\\ In the sequel $S$ is a finite set of $\Sigma_s$ which contains the set $S_{\infty}$ of archimedean places, $k_S$ the direct sum of the fields $k_s (s \in S)$ and $\mathcal{O}_S$ the ring of $S$-integers of $k$ (i.e. the ring of elements $x \in k$ such that $\vert x \vert_{s} \leq 1$ for $s \notin S )$.
  For $s$ non-archimedean, the valuation ring of the local field  $k_s$ is defined to be   
  $\mathcal{O}_{s} = \big\{ x \in k  \hh\hh |  \hh \hh \vert x \vert_{s} \leq 1 \big\}$. \\
%\textbf{Example.} If $k = \bbq$ and $p$ prime, for any $x\in \bbq $ of the form $x = p^n \frac{a}{b}$ with $a,b,n \in \bbz, b\neq 0$, the normalised absolute value at $p$ is given by $ \vert x \vert_{p} =  \vert p^n \frac{a}{b}\vert_{p} = p^{-n}$. The completion with respect to this absolute value is the field of $p$-adic numbers denoted by $\bbq_p$, its valuation ring is called the ring of $p$-adic integers and denoted by $\bbz_p$.\\
In all the statements of the article, without loss of generality one can replace $k$ by $\bbq$ but for sake of completeness we work with number fields.\\

\noindent Let $(Q,L)$ be a pair consisting of one quadratic form and one nonzero linear form on $k_{S}^{n}$. Equivalently, $(Q,L)$ can be viewed as a family $(Q_s,L_s) (s \in S)$, where $Q_s$ is a quadratic form on $k_{s}^{n}$ and $L_s$ a nonzero linear form on $k_{s}^{n}$. The form $Q$ is nondegenerate if and only each $Q_s$ is nondegenerate. We say that $Q$ is isotropic if each $Q_s$ is so, i.e. if there exists for every $s \in S$ an element $x_s \in k_{s}^{n} -\{0\}$ such that $Q_s(x_s) = 0$, in particular if $s$ is a real place an isotropic form is also said to be indefinite. For any quadratic form $Q$, we denote by rad$(Q)$ (resp. $c(Q)$) the radical (resp. the isotropy cone) of $Q$, by definition $Q$ is nondegenerate (resp. isotropic)  if and only if  rad$(Q) \neq 0$ (resp. c$(Q) \neq 0$).
The form $Q$ is said to be rational (over $k$) if there exists a quadratic form $Q_o$ on $k^n$ and a unit $c$ of $k_S$ such that $Q = c.Q_0$, and irrational otherwise. For any  $s \in S$ let $K_s$ denote an algebraic closure of $k_s$. If $G$ is a locally compact group, $G^{\circ}$ denotes the connected component of the identity in $G$.\\

\subsection{Main result} Let  be given a pair $F =(Q_s,L_s)_{s \in S}$ on $k_{S}^n$ and let $(a,b) \in k_{S}^{2}$.  We are interested in finding sufficient conditions which guarantees the existence of nontrivial $S$-integral solutions  $ x \in {\mathcal{O}}_{S}^{n}$ of the following simultaneous diophantine problem

\medskip

 \begin{center}
$(\mathrm{A_S})$  For any $ \varepsilon>0$,  $ |Q_{s}(x) - a_{s}|_s <  \varepsilon $  and $ |L_{s}(x) - b_{s}|_{s} <  \varepsilon $  for each $s \in S$. \end{center} 

\medskip 

\noindent Obviously as in the real case,  we need to find sufficient conditions on $F$ so that the set $ F({\mathcal{O}}_{S}^{n})$ would be dense in $k_{S}^{2}$. One have to be careful since the condition $(A_S)$ is not equivalent to density contrarily to real pairs (see \cite{BP}, \S6 and our \S\ref{sec end}).\\  

\noindent Our main result gives the required conditions for assertion $(A_S)$ to hold when $(a,b)= (0,0)$. In other words, we give sufficient conditions which implies that $ F({\mathcal{O}}_{S}^{n})$ is not discrete around the origin in  $k_{S}^{2}$.
 It may be seen as a weak $S$-arithmetic version of Theorem~\ref{thm},
\medskip

 \begin{thm}~\label{main thm} Let $Q = (Q_s)_{s \in S}$  be a quadratic form on $k_{S}^{n}$ and $L = (L_{s})_{s \in S}$ be a linear form on $k_{S}^{n}$ with $n \geq 4$ and $L_s \neq 0$ for all $s\in S$. Suppose that the pair $F=(Q,L)$ satisfies the following conditions, 
 
  \begin{enumerate}
\item $Q$ is nondegenerate.
\item $Q_{\vert L = 0}$ is nondegenerate and isotropic.
\item For each $s\in S$ the forms $ \alpha_s Q_s+ \beta_s L_{s}^{2}$ are irrational given any $\alpha_s,\beta_s$ in $k_{s}$ with $( \alpha_s,\beta_s) \neq (0,0)$.
\end{enumerate}

Then for any $ \varepsilon>0$,  there exists $x \in {\mathcal{O}}_{S}^{n}-\{0\}$ such that  
  \begin{center}
 $ |Q_{s}(x)|_s <  \varepsilon $  and $ |L_{s}(x)|_{s} <  \varepsilon $  for each $s \in S$.
 \end{center} 

\end{thm}

\medskip

\subsection{Remarks.}$(1)$ The proof of  Theorem \ref{main thm} reduces to dimension $4$, (see $\S ~\ref{sec reduc}$) this reduction is necessary since the proof relies essentially on classification of intermediate subgroups and requires low dimension\footnote{The classification of intermediate Lie subgroups in dimension greater than five becomes rapidly unfeasible when the dimension increases.}. The key is the use of the weak approximation in $k_S$ follows in the same lines (\cite{BP}, Proposition 1.3).

\medskip
\noindent  $(2)$ The proof of Theorem \ref{main thm} relies on Ratner's Theorem which gives a precise description of the closure orbits of lattices under the action of a connected Lie group generated by its unipotent one parameter subgroups. We need to apply an $S$-adic version of Ratner's theorem in order to find an integral solution simultaneously at all places. We treat first the case when $S=S_{\infty}$, by restriction of scalar we can use results of \cite{G1} to elucidate the structure of the intermediate subgroups. This is exactly where we need to work in dimension $4$, indeed the proof relies on the classifications of semisimple Lie algebras in $\mathfrak{sl}_4$ which contains the Lie algebra of the stabilizer.  For a general finite set of places $S$ containing both archimedean and nonarchimedean places, the use of strong approximation for number fields suffices to complete the picture.\\

\medskip

\noindent $(3)$ For Theorem \ref{main thm} even if we assume that $ \alpha Q + \beta L^{2}$ is irrational, it can be possible that the pencil form $ \alpha_s Q_s + \beta_s L_{s}^{2}$ is rational for some place $s$, in this situation it is not possible to apply Ratner's theorem. It can be possible that the result is still true in this situation but there are serious obstacles to (see \S~\ref{sec end}). 

\medskip
\noindent $(5)$ Unfortunately we are not able to show the density of $ F({\mathcal{O}}_{S}^{n})$ under the conditions of theorem~\ref{main thm} with our method. We are also even unable to show that $ |Q_{s}(x)|_s $  and $ |L_{s}(x)|_{s}$ are both nonzero for any $s \in S$  and $x \in {\mathcal{O}}_{S}^{n}$ as in the conclusion of Theorem~\ref{main thm}. We discuss those issues in  $\S~\ref{sec end}$.\\

\noindent $(6)$ In the real case, one can hope to relax condition $(2)$ by only asking  $ \alpha Q + \beta L^{2}$  to be isotropic as it is conjectured by Gorodnik (see \S~\ref{sec end}, Conjecture $\ref{conj} $). The major issue is that reduction to lower dimension fails to hold.

\subsubsection*{Acknoledgements.} This article is  part of the PhD thesis of the author under the supervision of A.Ghosh. I thank my supervisor for suggesting me this problem, and for many fruitful discussions.

\section{Weak and strong approximation, reduction to dimension 4}~\label{sec reduc}
\subsection{Weak approximation in number fields and Grassmannian varieties}
\medskip

\noindent Number fields satisfy a nice \textit{local-global principle} called the weak approximation which can be seen as a refinement of the Chinese remainder theorem. 
\begin{thm}[Weak approximation in number fileds] Let $S$ be a finite set of  $\Sigma_k$. Let given $\alpha_s \in k_s$ for each $s \in S$. Then there exists an $\alpha \in k$ which is arbitrarily close to $\alpha_s$ for all $s \in S$ with respect to the $s$-adic topology. 
\end{thm}
\noindent \textbf{Proof.} (See e.g. \cite{L94}, Theorem1, p.35)\\

One can reformulate this theorem as follows: the diagonal embedding $ k \hookrightarrow \prod_{s \in S} k_s$ is dense, the product being equipped with the product of the $s$-adic topologies. 

\begin{definition}[Weak approximation in algebraic varieties]  Let $X$ be an algebraic variety defined over $k$, then $X$ is said to satisfies weak approximation property with respect to $S$ if the diagonal embedding $X(k)  \hookrightarrow \prod_{s \in S} X(k_s)$ is dense for the $S$-adic topology.

\end{definition}
To prove reduction we need to introduce a useful class of algebraic varieties which satisfies weak approximation,
\begin{definition} Let $V$ be a $k$-vector space of dimension $n\geq 1$ and for each $1\leq m\leq n$ let us define the set
\begin{center}
$ \mathcal{G}_{m}(V) = \Big\{k$-vector subspaces $W$  $ \subset V$ with $\dim W = m \Big\}$.
\end{center}
This is an algebraic variety defined over $k$ called the Grassmannian variety, if $V=k^n$ we just denote it  as $ \mathcal{G}_{n,m}(k)$. 
\end{definition}

\begin{prop} Let be given two integers $1 \leq m \leq n$, then the Grassmannian variety  $\mathcal{G}_{n,m}(k)$  satisfies weak approximation with respect to $S$, that is,
\begin{center}
$  \mathcal{G}_{n,m}(k)  \hookrightarrow \prod_{s \in S}  \mathcal{G}_{n,m} (k_s)$ is dense. 

\end{center}
\end{prop}
\noindent \textit{Proof.} Let be given a family $(V_s)_{s\in S}$ of $k_s$-vector subspaces of dimension $m$ in $k_{s}^{n}$ for each $s\in S$. Each of these subspaces $V_s$ are determined by $m$ linearly independent vectors in $k_{s}^n$. For each of the $V_s$, the coefficients of these vectors in the standard basis of $k_{s}^n$ give rise to a $m\times n$-matrix $A_s$ with coefficients in $k_{s}$. By weak approximation property in $k_{S}^{nm}$ we obtain a matrix $B \in \mathcal{M}_{m,n}(k)$ such that for any $s\in S$, $B_s$ is arbitrarily close to $A_s$. Let $V^{^{\prime}}$ be the vector subspace generated by the $n$ columns of $B$, obviously $V\in \mathcal{G}_{n,m}(k)$ and $V^{\prime}_s$ is arbitrary close to $V_s$ for all $s\in S$.

 \subsection{Reduction of Theorem \ref{main thm} to the dimension 4}

 \begin{prop}
Let $ F = (Q,L)$ be a pair consisting of a quadratic form $Q$ and a nonzero linear form $L$ in $k_{S}^{n}$ $(n \geqslant 5)$ such that 
\begin{enumerate}
  \item $Q $ is nondegenerate
  \item $Q_{| L = 0}$ is isotropic
  \item  Any quadratic form $ \alpha_s Q_s + \beta_s L_{s}^{2}$  with  $\alpha_s,\beta_s$ in $k_{s}$ such that $( \alpha_s,\beta_s) \neq (0,0) $ for all $s \in S$  is irrational.  \end{enumerate}
Then there exists a $k$-rational subspace $V$ of $k^{n}$ of codimension $1$ such that $F_{|V_{S}}$ satisfies the conditions $(1)(2)(3)$, moreover $V$ can be chosen such that $Q_{| \{ L = 0 \} \cap V_S }$ is nondegenerate.
\end{prop}
\noindent \textbf{Proof.} When $s$ is an archimedean real place, it is proved in (\cite{G1}, Proposition 4) that there exists a subspace $V_{s}$ of $k^n$ of codimension 1 such that $F_{s |V_{s}}$ verifies conditions $(1)(2)(3)$. In the case of archimedean complex places and nonarchimedean places, one may replace the condition $Q_{s | L_{s} = 0 }$ of type (I) which only  valid for real places by equivalent condition that $Q_{s | L_{s} = 0 }$ is nondegenerate which is valid for all $s \in S$. Therefore there exists a subspace $V_{s}$ of $k^n$ of codimension 1 such that $F_{s |V_{s}}$ verifies conditions $(1)(2)(3)$, the proof of the latter existence of $V_s$ for non-archimedean places in $S$ is identical to the real places (see \cite{G1}, Proposition 4). Hence for any $s \in S $ we may find  $V_{s}$ a subspace of $k^n$ of codimension 1 so that the conditions $(1)(2)(3)$ are satisfied by $F_{s |V_{s}}$ and one can choose $V_s$ to be such that $Q_{s | \{ L_{s} = 0 \} \cap V_{s}}$  is nondegenerate. \\
Assume that $n \geqslant 5$. Let be given $s \in S$ and  $V_s$ a $k$-subspace of codimension 1 in $k_{s}^n$ such that the restriction of $Q_s$ on $V_s$ is non-degenerate and isotropic. Let us define $H_s :=  SO(Q_{s}) $ the $k_s$-algebraic subgroup of the orthogonal group, $H_s(K_s)$ is a connected Lie group over the algebraic closure $K_s$ of $k_s$. F

\noindent Following Borel and Prasad (\cite{BP} Proposition 1.3),   let us consider the $k_s$-action of  the group $H_s$ on the Grassmanian variety ${\gcal}_{n-1,n}$ of the hyperplanes over $k_s$.
Then we use the fact that the orbit $H_s(K_s)V$ is open in ${\gcal}_{n-1,n}(K_s)$  for the analytic topology.  The fact that the fibration $\pi$ of  the orbit of $V_s$ under  $H_s(k_s) $ by the isotropy group of $V_s$ has local $k_s$-cross-section $\sigma$ implies that $H_s(k_s) V_s$ is also open in ${\gcal}_{n-1,n}(k_s)$ for the analytic topology.

\begin{center}
\begin{tikzpicture}

\node (st) at (-1.5,3) {$Stab_G(V_s)  \longrightarrow$};
\node (S) at (0,3) {$H_s $};
\node (Y) at (3,3) {$H_s.V_s  $};
\node (pi) at (1.6,3.2) {$\pi$};
\node (sig) at (1.6,2) {$\sigma$};
\node (G) at (4.5,3) {$ \hookrightarrow {\gcal}_{n-1,n}$};
\tikzstyle{estun}=[->,dotted,very thick,>=latex]
\draw[estun] (Y) to[bend left] (S);
\draw[->] (S) to (Y);

%\draw[estun] (st) to[right] (S);
\end{tikzpicture}

\end{center}

%$ \mathrm{SO}(Q_{ s}) \mapsto {\gcal}_{n-1,n}(k_s)$ given by $g \mapsto g(V_s)$ is submersive\footnote{i.e. the induced map on tangent space is surjective.} 
%\end{center}
%in particular the image contains a neighbourhood of the image of any $g$ after the implicit function theorem (see e.g.\cite{ev}, Lemma 2). 
\noindent Moreover by weak approximation in $k_S$ we can find a rational subspace in $V^{\prime}$ of codimension 1 in $k^n$ such that $V^{\prime} \otimes_k k_s$ is arbitrarily close to $V_s$ for all $s \in S$, in particular they belong to the same open orbit under $H_s$. We have established that  $F_{s |V_{s}}$ satisfies conditions $(1)$ and $(2)$,  it is equivalent to say that 
\begin{center}
rad$(Q_s ) \cap V_s = \{0\}$  and $c(Q_{ s|L_{s} =0}) \cap V_s \neq \{0\}$.         $(\ast)$ \end{center}

\noindent The condition $(2)$ remains true if we replace $V_s$ by any subspace sufficiently close to $V_s$. Since the subspace rad$(Q_s)$  is invariant under the action of the orthogonal group $\mathrm{SO}(Q_{ s})$, the condition $(1)$ above is verified by any element of ${\gcal}_{n-1,n}(k_s)$ which lies in the orbit of $V_s$ under  $\mathrm{SO}(Q_{ s})$. In particular, $V^{\prime} \otimes_k k_s$  satisfies $(\ast)$ for each $s \in S$. Hence we obtain a $k$-rational subspace $V^{\prime}$ of $k^n$ such that $F_{|V^{\prime}_{S}}$ satisfies the conditions $(1)(2)$. It remains to find such $V$ such that in addition $F_{|V_{S}}$ satisfies condition $(3)$. 
\noindent  Let us put
\begin{center}
$\V = \Big\{V \in  {\gcal}_{n-1,n}(k) \hh | \hh   \hh  \hh F_{ |V_{S}} \hh  \hh  \mathrm{satisfies} \hh \hh \mathrm{conditions} \hh$(1)\hh(2)$ \Big\}  $.

\end{center}
It is nonempty because it contains $V^{\prime}$. Suppose there exists no $V$ in $\V$ for which $F_{|V_{S}}$ satisfies condition (3), that is to say that for any $V\in \V$, it should exists some $ s \in S$ and some $(\alpha_s, \beta_s) \in k_{s}^2 -\{(0,0)\}$, such that the quadratic form $\alpha_s Q_{s}(x)+ \beta_{s} L_s(x)^2_{|V_{s}}$ is rational.  Let us consider the regular map $ f: k_{s}^n \longrightarrow k_s$ given by $$f: x \mapsto \alpha_s Q_{s}(x)+ \beta_{s} L_s(x)^2.$$

\noindent Clearly $f$ is a polynomial function on $K_{s}^n$. For each $V \in \V$ we have  $f(V(k)) \subset k$ and the Zariski density of  $V(k)$ in  $\overline{k_s}^n$ implies that  $f$ is defined over $k$.  In other words,  $\alpha_s Q_{s}(x)+ \beta_{s} L_s(x)^2$ is rational over $k$, contradiction. Hence there exists $V \in \V $ such that  $F_{|V_{S}}$ satisfies condition $(3)$. 
\begin{cor}
It suffices to prove Theorem 2.1 for $n = 4$.
\end{cor}
\noindent \textbf{Proof.} It follows from the proposition by descending induction on $n$.

\subsection{Adeles and strong approximation for number fields} The set of adeles  $\mathbb{A}$ of $k$ is the subset of the direct product $\prod_{s \in \Sigma_k} k_s$ consisting of those $x =(x_s)$ such that $ x \in \mathcal{O}_s$ for almost all $s \in \Sigma_k$.  The set of adeles $\mathbb{A}$ is a locally topological ring with respect to the adele topology given by the base of open sets of the form $ \prod_{s \in S} U_s   \times  \prod_{s \notin S} \mathcal{O}_s$ where $S \subset \Sigma_k$ is finite with $ S \supset S_{\infty}$ and $U_s$ are open subsets of $k_s$ for each $s \in S$. For any subset $S \subset \Sigma_k$ finite with $ S \supset S_{\infty}$, the ring of $S$-integral adeles is defined by: 
 \begin{center}
 $\mathbb{A}(S)  =   \prod_{s \in S} k_s   \times  \prod_{s \notin S}  \mathcal{O}_s$, thus we can see that $\displaystyle \mathbb{A}=  \bigcup_{S \supset S_{\infty}} \mathbb{A}(S)$.
\end{center}
We define also $\mathbb{A}_S$ to be the image of $\mathbb{A}$ onto $\prod_{s \notin S} k_s$, clearly $\mathbb{A} = k_S \times \mathbb{A}_S$.
\begin{thm}[Strong approximation] If $S \neq \emptyset$ the image of $k$ under the diagonal embedding is dense in $\mathbb{A}_S$.
\end{thm}

\section{Stabilizers of pairs $(Q,L)$}~\label{sec stab}

\noindent For  each $s \in S$ let us define $G_{s} = \mathrm{SL}_{4}(k_s)$,  $G_{S} =\prod_{s \in S} \mathrm{SL}_{4}(k_s) = \mathrm{SL}_{4}(k_S)$. 
Let $F=(Q,L)$ be a pair on $k_{S}^{4}$ satisfying  the conditions $(1)(2)(3)$ of Theorem ~\ref{main thm}.

\noindent For every $s \in S$ we realize $Q_s$ on a four-dimensional quadratic vector space $(W_s,Q_s)$ over $k_{s}$  equipped with the standard basis $ \Bc = \{e_1,\cdots,e_4\}$. For each $s \in S$, let us define $H_s$ the stabilizer of the pair $F_s$ under the action of $G_s$, in other words 
\begin{center}
$H_s = \Big\{ g \in G_s  \hh \hh | \hh \hh   Q_s \circ g = Q_s \hh, \hh  L_s \circ g = L_s \Big\}$.
\end{center}
Equivalently one can write $H_s = \Big\{ g \in \mathrm{SO}(Q_s)  \hh \hh | \hh  \hh \hh  L_s \circ g = L_s \Big\}$, clearly it is a linear algebraic group defined over $k_s$. Also let us define $V_s =\{L_{s}=0\}$, it is an hyperplane of $W_s$ which induces a quadratic isotropic subspace $(V_{s},Q_{ s| V_{s} })$ of dimension $3$ in $W_s$. We have two cases following $(V_{s},Q_{ s| V_{s} })$ is nondegenerate or not. If $s$ is a real place the first case corresponds to pairs of type (I) in the terminology of \cite{G1}. \\

\begin{lem}\label{lemma} Let be given a pair $(Q,L)$ satisfying the conditions of Theorem \ref{main thm} in dimension 4. Then the stabilizer of $(Q,L)$ under the action of $G$ is of the form (up to conjugation) 
\begin{center}
$H = \bigg\{  \left(\begin{array}{c|c} A & 0 \\\hline 0 & 1\end{array}\right)  \hh \hh \hh  | \hh \hh \hh A \in \mathrm{SO}(Q_{ | L=0 })   \bigg\} \subseteq \hh \hh\mathrm{SL}_{4}(\overline{k_{s}})$.\end{center}
In particular, $H$ is semisimple. Moreover, any quadratic form $\widetilde{Q}$ which is $H$-invariant is of the form $\alpha Q +\beta L^2$ for some $\alpha,\beta \in k_S$ not both zero.
\end{lem}
\noindent \textit{Proof.} Since $(V_{s},Q_{ s| V_{s} })$ is nondegenerate, one can write the the following decomposition $W_s = V_s \oplus V_{s}^{\bot}$ where $V_{s}^{\bot} $ the orthogonal complement w.r.t. $Q$. Since $\dim V_{s}^{\bot}=1 $ there exists some nonzero $u$ vector of $W_s$ such that $V_{s}^{\bot} = <u>$ with $L_s(u) \neq 0$. Moreover by definition $L_s$ is $H_s$-invariant so $V_s$ is $H_s$-invariant. Moreover any element of $h\in H_s$ is in particular an element of $SO(Q_s)$, that is, $h^{T} = h$ hence $ V_{s}^{\bot}=<u>$ is also $H_s$-invariant. Then for any $h\in H_s$,  the restriction $h_{|V_s} $ induces an automorphism of $V_s$ and $h u =u$. Let us put $w_4=u$, and complete with a basis of $V_s$ $\Big\{w_1,w_2,w_3\Big\}$ to obtain the following matrix representation of $H_s$ up to an $k_s$-isomorphism of $W_s$,

\medskip

% Since $V_{s}^{\bot} $ is clearly $H_s$-invariant, therefore $H_s$ acts by $x \mapsto \lambda x$ on the line $<v>$. The linear form $L_{s| V_{s}^{\bot}}$ is nonzero and is $H_s$-invariant then $H_s$ acts trivially on the line $V_{s}^{\bot}$. Let us define $w_4 = v$ and complete to obtain a basis ${\Bc}^{\prime} = \{w_1,\cdots,w_4\}$ of $W_s$ where $<w_1,w_2>$ is an hyperbolic plane in $V_s$. 
%Hence   we obtain for any $x \in W_s$, with respect to the basis ${\Bc}^{\prime} = \{w_1,\cdots,w_4\}$  that
 %\begin{center}
% $Q_{s}(x) =  x_{1}x_{2}  + a_{3}x_{3}^2 + a_{3}x_{4}^2 $ and $L_{s}(x) = x_4 $ with  $a_3, a_4 \in k_{s}^{\ast}$.
%\end{center}
\begin{center}
$H_s \simeq  \bigg\{  \left(\begin{array}{c|c} A & 0 \\\hline 0 & 1\end{array}\right)  \hh \hh \hh  | \hh \hh \hh A \in \mathrm{SO}(Q_{s | V_{s} })   \bigg\} \subseteq \hh \hh\mathrm{SL}_{4}(\overline{k_{s}})$.
\end{center}
It is well-known that the orthogonal group of a nondegenerate quadratic form is a semisimple group. The last statement is the Lemma 9 in \cite{G1}.

\subsection*{$S$-adic products}

\noindent Now let  $F=(Q_s,L_s)_{s \in S}$ be a pair satisfying the conditions of the main theorem.
 Let $\mathcal{H}_s$ be the algebraic group defined over $k_s$ such that $\mathcal{H}_s(k_s) = H_s$. Define $H_{s}^{+}$ to be the subgroup of $H_s$ generated by its one-dimensional unipotent subgroups. 
 \noindent Let us put \begin{center}
$H_S = \prod_{s \in S} H_s$ and $H_{S}^{+} = \prod_{s \in S} H_{s}^{+}$.
\end{center}

\noindent Therefore $H_S$ is an algebraic subgroup of $\mathrm{SL}_{4}(k_S)$ which leaves invariant the pair $F = (Q,L)$ with respect to the $S$-basis $\mathcal{B}^{\prime} = \bigg\{w_1,w_2,w_3,w_4\bigg\}$ as in the previous lemma.  \\ In other words, we have

\begin{center}
$H_S \simeq  \bigg\{  \left(\begin{array}{c|c} A & 0 \\\hline 0 & 1\end{array}\right)  \hh \hh \hh  | \hh \hh \hh A \in \mathrm{SO}(Q_{ | V_S })   \bigg\} \subseteq \hh \hh\mathrm{SL}_{4}(\overline{k_{S}})$.
\end{center}
and

\begin{center}
$H^{+}_{S} \simeq  \bigg\{  \left(\begin{array}{c|c} A & 0 \\\hline 0 & 1\end{array}\right)  \hh \hh \hh  | \hh \hh \hh A \in \mathrm{SO}(Q_{ | V_S })^{+}   \bigg\} \subseteq \hh \hh\mathrm{SL}_{4}(\overline{k_{S}})$.
\end{center}

%Let be given a quadratic form $\widetilde{Q}$ in $k_{S_1}^{4}$ and assume it to be invariant under $H_1$. Let $s \in S_1$, one can write 
%$\widetilde{Q}_s(x) = \widetilde{Q}_{s| L = 0}(x) + \widetilde{L}_s(x_1,x_2,x_3)x_4 + \widetilde{Q}_s(w_4) x_{4}^{2}$ for some linear form $\widetilde{L_s}$ in $k_{s}^{3}$. Since $\widetilde{Q}_s$ is $H_s$-invariant, $ \widetilde{Q}_{s| L = 0}$ and $ \widetilde{ L}_s$ are both invariant under the orthogonal group $\mathrm{SO}(Q_{s | L_{s}=0 })$. Hence $ \widetilde{Q}_{s| L = 0}$ is scalar multiple of $Q_{s | L_s = 0}$ and $\widetilde{L}_s = 0$ thus we obtain the analog of (\cite{G1} lemma ?)
%\begin{lem}
%Let  $\widetilde{Q}$ be a quadratic form in $k_{S_1}^{4}$ which is invariant under $H_1$. Then $\widetilde{Q} = \alpha Q + \beta L^{2}$ for some $\alpha,\beta \in k_{S_1}$. 

%\end{lem} 

% \subsubsection*{Intermediate subgroups of  $G_{S}$} 

\section{ Topological rigidity in  $S$-adic homogeneous spaces}~\label{sec uni} 
  \subsection{Ratner's topogical rigidity theorem for unipotent groups actions}
 \noindent  Let $G_{S} = \mathrm{SL}_{4}(k_S)$ and let $\Gamma_S$ be the $S$-arithmetic subgroup of $G_S$ given by $\Gamma_S =    \mathrm{SL}_{4}({\mathcal{O}}_S)$.  The ring $\mathcal{O}_S$  is a lattice in $k_S$. Let us define $\Omega_S$ to be the quotient space given by $G_S/\Gamma_S$. It is the space of free of ${\mathcal{O}}_S$-submodules of $k_{S}^{4}$ of maximal rank and determinant one. Then $\Omega_S$ is the homogeneous space of unimodular lattices of ${\mathcal{O}}_{S}^{4}$, by lattice we mean a discrete subgroup of $G_S$ of finite covolume. For every $s \in S$, let $\mathcal{U}_s$ be a unipotent $k_s$-algebraic subgroup of $\mathrm{SL}_{4/k_s}$ and denote by $\mathcal{U} = \prod_{s \in S} \mathcal{U}_{s}(k_s)$ the associated unipotent subgroup of $G_S$.
 
\noindent We are interested in the left action of $\mathcal{U}$ on the homogeneous space $\Omega_S$ and more particularly with the closure of such orbits. If $x \in \Omega_S$ it turns out that the closure of the orbit $\mathcal{U}x$ is also an orbit of $x$. The following result is the generalisation of Ratner's orbit closure theorem for $S$-products proven independently by Margulis-Tomanov and Ratner (see \cite{MT}, \cite{R}). 
 
 \medskip
\begin{thm}[Ratner's Theorem for $S$-adic groups]~\label{MT} Assume that  $ \mathcal{U}$ is generated by its one-dimensional unipotent subgroups. Then for any $x \in \Omega_S$, there exists a closed subgroup $M = M(x) \subset G_S $ containing $\mathcal{U}$ such that the closure of the orbit $\mathcal{U}x$ coincides with $Mx$ and $Mx$ admits $M$-invariant probability measure.
\end{thm}

\noindent In the real case, this result was conjectured by Raghunathan who stated it with $G=SL(3,\bbr)$, $U=SO(2,1)$ and $\Gamma = SL(3,\mathbb{Z})$. He noticed that the proof of this Conjecture gives a  the Oppenheim conjecture.  Despite appearences, the proof of this theorem is measure theoretic and consists to classify ergodic measures under the action of the unipotent flow on the homogeneous space $\Omega$.  In both proofs made by Margulis-Tomanov (\cite{MT}) and Ratner (\cite{R}) the notion of entropy plays a central role for the measure classification.

\subsection{The structure of intermediate subgroups in the archimedean case.}

The application of Ratner's Theorem \ref{MT} above  gives a nice description of the orbit closure under the action of a subgroup $H^{+}$ generated by unipotent one parameter subgroups in $G_S$. Considering a unimodular lattice $x \in \Omega_S$, we get that 
$$\overline{H^{+}x} = L. x$$
for some closed connected subgroup $L$ (depending on $x$) such that $H^{+} \subset  L \subset G$.
 Indeed, the solution of the Oppenheim conjecture relies essentially in proving that the orbit closure of any element $G/\Gamma$ under the action of $SO(2,1)^{\circ}$ is either closed or dense in $G/\Gamma$.
\medskip
\noindent The following theorem due to N. Shah gives  extra information about the structure of the intermediate subgroups arising from Ratner's orbit closure theorem above, 
\begin{thm}[\cite{Sh}, Prop. 3.2]~\label{shah} Let $G = \mathcal{G}(\mathbb{R})^{\circ}$ with $\mathcal{G}$ an algebraic subgroup of $SL(n, \mathbb{C})$ defined over $\bbq$ and $\Gamma = G \cap SL(n,\mathbb{Z})$. Let $H$ be a subgroup of $G(\bbr)$ generated by its unipotent one parameter subgroups and assume that 
\begin{center}
$\overline{H\Gamma} = P\Gamma$ where $P$ is a closed connected subgroup of $G(\bbr)$\end{center}
 such that $P \cap \Gamma$ has finite covolume in $P$. Then $P = \widetilde{P}(\bbr)^{\circ}$ where $\widetilde{P}$ is the smallest $\bbq$-subgroup of $G$ whose group of real points contains $H$.
\end{thm}
\noindent By combining Ratner's theorem and previous proposition, we obtain immediately the following corollary stated in the set-up of the previous proposition.
\begin{cor}[ \cite{BP}, Proposition 7.2]\label{corP}  Let $g \in G$ such that $x=g.0$ where $0$ is the coset of $\Gamma$ in $G/\Gamma$. Then $g^{-1}Lg = \widetilde{P}(\mathbb{R})$ where $\widetilde{P}$ is the smallest subgroup of real points which contains $g^{-1}Hg$.

\end{cor}
\noindent To use the results of this section for our purposes, we need to reduces to the case when $k=\bbq$. This can be done by using the functor of restriction of scalars for algebraic groups.
\begin{prop}
Let $k$ be a number field.  Given any algebraic group $G \subset GL_n(\overline{k})$ defined over $k$, there exists an algebraic group $G^{\prime}$ defined over $\bbq$ such that $G^{\prime}(\bbq) \simeq G(k)$.
\end{prop}
\noindent \textit{ Proof.}  See e.g.  \cite{PR}, $\S 2.1.2$ for a general construction for finite separable extensions.
\begin{definition}[Restriction of scalars]

We denote $G^{\prime}$ by $\mathcal{R}_{k/\bbq} (G)$ and it is called the algebraic group associated to $G$ obtained by restriction of scalars from $k$ to $\bbq$.
\end{definition}
\noindent The operation $\mathcal{R}_{k/\bbq}$ defines a functor from the category of $k$-groups to the category of $\bbq$-groups. This functor has a nice arithmetic property regarding the set of integral points, given any algebraic group defined over $k$, we have\footnote{This result will not be used in the sequel.} 
$$ \mathcal{R}_{k/\bbq}(G)(\mathbb{Z}) \simeq G(\mathcal{O}_k)$$
Assume $G=SL_{n\mid k}$ view as a algebraic group defined over $k$ and if  all the places are archimedean i.e. $S=S_{\infty}$. The fact that $k_{S} = k \otimes_{\bbq}\bbr$  implies $G_{S} =  \mathcal{R}_{k/\bbq} (G) (\bbr)$ and the intermediate subgroups have a still a interesting structure by means of restriction of scalars

\begin{prop}[\cite{BP}, Prop. $7.3$]\label{res} Let $H_s$ be a closed subgroup of $G_s = SL_n(k_s)$ for each $s \in S_{\infty}$ and $H$ the product of the $H_s$. Then the smallest $\bbq$-algebraic subgroup $\mathcal{L}$ of $\mathcal{G}$ whose group of real points contains $H$ is of the form $\mathcal{L} = \mathcal{R}_{k/\bbq} \mathcal{L}^{\prime}$, where 
$ \mathcal{L}^{\prime}$ is a connected $k$-subgroup of $\mathcal{G}$.
\end{prop}

\section{Proof of the Theorem~\ref{main thm}} 
\noindent Let $F= (Q,L)$ be a pair in $k_{S}^{n}$ which satisfies the conditions of Theorem \ref{main thm}. After $\S~\ref{sec reduc}$, we know that it suffices to show it for $n = 4$. 
By condition $(3)$ all the forms $ \alpha_s Q_s + \beta_s L_{s}^{2}$ are irrational for each  $\alpha_s,\beta_s$ in $k_{s}$ such that $( \alpha_s,\beta_s) \neq (0,0)$ for any $s \in S$. Let $g \in G_S$ be the matrix of the basis $\mathcal{B}^{\prime}_{S}$ in the standard basis of $k^{4}_{S}$. By definition $g^{-1}H_S g$ leaves invariant the pair $F=(Q_s,L_s)_{s \in S}$,  and $H_{S}^{+}$ is generated by one-dimensional unipotent subgroups. We consider $\Gamma_S$ as an element of the homogeneous space $\Omega_S$. By applying Ratner's Theorem \ref{MT}, one obtains
 \begin{equation}\label{ratner}
 \overline{g^{-1}H_{S}^{+}g \Gamma_S} = P \Gamma_S
\end{equation}
where $P$ is a closed subgroup of $G_S$ which contains $g^{-1}H_{S}^{+} g$.
 %Moreover  there is a Levi decomposition $\widetilde{M} = M_0 R_{u}(\widetilde{M})$ where  $R_{u}(\widetilde{M})$ the unipotent radical of $\widetilde{M}$ and $M_0$ a connected semisimple algebraic subgroup of $M$.

\medskip
\noindent Assume first that $S = S_{\infty}$, thus $\mathcal{O}_{S}^{4} = \mathcal{O}^4$  and we simply write $k_{\infty}$, $H_{\infty}$ and $G_{\infty}$ respectively for $k_{S_{\infty}}$,  $H_{S_{\infty}}$ and $G_{S_{\infty}}$. One notes also that $H_{\infty}^{+}$ is nothing else than the component of the identity $H_{\infty}^{\circ}$. Using equality (\ref{ratner}) one deduces that the set  $F({\mathcal{O}}^4)$ is dense in $k_{\infty}^2$. Indeed, we are going to adapt the proof of (\cite{G1}, Proposition 10) to the $S_{\infty}$-products, as follows\footnote{For more precisions, the reader is invited to read the original proof which is similar.}. We first reduce the ground field from $k$ to the field of rational numbers. To achieve this we realise $G_{\infty}$ as the group of real points of an algebraic group $\mathcal{G}$ defined over $\bbq$. In view of Proposition \ref{res} this is given explicitly by taking $\mathcal{G}= R_{k/\bbq} \SL_{4}$ where $R_{k/\bbq}$ is the functor restriction of scalars of the field extension $k/\bbq$ and where $\SL_{4}$ is regarded as the usual algebraic group over $k$. In other words, $G_{\infty} = \mathcal{G}(\bbr)$ with $\mathcal{G}$ an algebraic group defined over $\bbq$. Now let us precise the structure of $P$. From Corollary \ref{corP} above, we infer that there exists an algebraic group $\widetilde{P}$ defined over $\bbq$ which is the smallest $\bbq$-algebraic group whose group of real points $\widetilde{P}$ contains $g^{-1}H_{\infty}^{\circ}g$. In the other hand, Proposition \ref{shah} implies that $P= \widetilde{P}(\bbr)^{\circ}$ and the unipotent radical $U$ of $\widetilde{P}$ is also defined over $\bbq$. Thus equality (\ref{ratner}) may be read as
\begin{equation}\label{rshah}
 \overline{g^{-1}H_{\infty}^{\circ}g \Gamma} = \widetilde{P}(\bbr)^{\circ} \Gamma.
 \end{equation}
\begin{lem}\label{archimedean} For each $s\in S_{\infty}$, let $P_s$ be the intersection of $P$ with $G_s$. If $P_s$ acts irreducibly on $\bbc^4$, then $P_s = G_s$. Otherwise,  $P_s = M_{s}U$ where 
\begin{center}
$M_s =  ug_{s}^{-1}\left(\begin{array}{c|c} \SL_3 & 0 \\\hline 0 & 1\end{array}\right)g_{s}u^{-1}$ for some $u \in U_s$.

\end{center} 

\end{lem}

\noindent \textit{Proof the Lemma.} This result is the core of the proof of Proposition 10 in \cite{G1} for which we recall the outlines. If $P_s$ acts irreducibly on $\bbc^4$, then $P_s$ is semisimple and the classification of irreducible semisimple Lie groups in $\SL_4$ implies that $P_s$ is equal either to $G_s$ or $SO(B_s)$ for some nondegenerate form $B_s$ (Proposition 7 and Lemma 8, \cite{G1}). Such form $B_s$ being $H_s$-invariant is necessarily of the form $ \alpha_s Q_s + \beta_s L_{s}^{2}$ for some $( \alpha_s,\beta_s) \neq (0,0)$ (Lemma \ref{lemma}). As seen before $P_s$ is defined over $\bbq$, so that $B_s$ is forced to be rational which is a contradiction. Hence $P_s = G_s$. For the second assertion, we consider the induced action of $P_s$ on the space $\mathcal{L}$ of linear forms in $\bbc^4$, it is reducible by hypothesis.  There are only two $P_s$-invariant subspaces in $\mathcal{L}$, namely $\mathcal{L}_1 = \langle L_1, L_2, L_3\rangle$ and  $\mathcal{L}_2 = \langle L_4\rangle$ where $L_i(x) =(gx)_i$ for $i=1,\ldots,4$, note that $L_4 = L$. Since $P_s$ is defined over $\bbq$, one infers that $M$ is semisimple thus admitting a Levi decomposition 
\begin{center}
$P_s = M_{s} U_s$  
\end{center}
where $M_s$ and $U_s$ are respectively a Levi subgroup and the unipotent radical of $P_s$. The Levi subgroup $M_s$ is defined over $\bbq$ since $P_s$ is. Also as seen above, $U$ is defined over $\bbq$ and Malcev's theorem ensures that the Levi subgroups are unique up to conjugacy (e.g. see \S4.3 \cite{OV}), in particular 
\begin{equation}\label{strict}
g_{s}^{-1} H_{s}^{\circ}g_s \subseteq u^{-1}M_{s}u \hh \hh \mathrm{for \hh some }\hh \hh  u\in U_s.
\end{equation}
Moreover this inclusion is strict because $H_s$ is not defined over $k$.  The latter fact and the maximality of $\mathrm{SO}(Q_{\mid L=0})$ in $ \SL_3$ ($Q_{\mid L=0}$ is isotropic) gives the equality
\begin{center}
$M_s =  ug_{s}^{-1}\left(\begin{array}{c|c} \SL_3 & 0 \\\hline 0 & 1\end{array}\right)g_{s}u^{-1}$.
\end{center}
This achieves the proof of the Lemma.\\

\noindent Let us define the subgroup 
\begin{center}
$M^{\prime}_s :=u^{-1}M_{s}u =  g_{s}^{-1}\left(\begin{array}{c|c} \SL_3 & 0 \\\hline 0 & 1\end{array}\right)g_{s}$.
\end{center}
By the previous Lemma \ref{archimedean}, one can rephrase equality (\ref{rshah}) in the following way

\begin{equation}\label{density}
 \overline{g^{-1}H_{\infty}^{\circ}g \Gamma} = M^{\prime}(\bbr)^{\circ} U(\bbr)\Gamma.
 \end{equation}

\noindent Now let be given $(a,b) \in k_{\infty}^2$ and let us choose $x\in \mathcal{O}^4 - \langle g^{-1}e_4\rangle$. It is not difficult to see that there exists $m \in M^{\prime}(\bbr)^{\circ}$ and $u\in U(\bbr)$  such that
\begin{center}
$F(mux) =(Q(mux),L(mux)) =(a,b)$.
\end{center}
Using density in (\ref{density}), we infer that there exists $h_n \in g^{-1}H_{\infty}^{\circ}g$ and $\gamma_n \in \Gamma$ such that 
\begin{center}
$h_{n}\gamma_{n} \rightarrow mu$ as $n \rightarrow \infty$.
\end{center}
We conclude that \begin{center}
$F(\gamma_{n}x)= F(h_{n}\gamma_{n}x) \rightarrow F(umx) =(a,b)$ as $n \rightarrow \infty$.
\end{center}
In other words,  $F(\mathcal{O}^4)$ is dense in $k_{\infty}^2$ and in particular this proves Case 1 when $S = S_{\infty}$.\\

\noindent Now let us assume $S \neq S_{\infty}$ and let be given $s \in S_f$. The set ${\mathcal{O}}^4$ is bounded in $k_{s}^{4}$, thus for any neighbourhood $U$ of the origin in $k_{s}^{4}$ one can find an integer $a_s \in {\mathcal{O}}_s $ such that $a_s .{\mathcal{O}}^4 \subset U$.  In other words, given any $\varepsilon > 0$ one can find $a_s \in  {\mathcal{O}}_s$ such that :
\begin{center}
 $ |Q_{s}(a_s x)|_s \leqslant \varepsilon$ and $ |L_s( a_s x)|_s \leqslant \varepsilon$  for all $x \in {\mathcal{O}}^4$.
\end{center}
Thus for each $s \in S_f$, we can associate an integer $a_s \in  {\mathcal{O}}_s$ satisfying the previous inequalities. By strong approximation one can find $a \in  {\mathcal{O}}$ such that $|a|_s = |a_s|_s$ for all $s \in S_f$. Put $\| a \|_{\infty} = \max_{s \in S_{\infty}} |a|_s$, by the previous case we can find $x \in {\mathcal{O}}^4$ such that: 
\begin{center}
 $ |Q_{s}(x)|_s \leqslant \varepsilon/ \| a \|_{\infty}^{2}$ and $ |L_s( x)|_s \leqslant \varepsilon/ \| a \|_{\infty} $  for all $s \in S_{\infty}$.
\end{center}
We immediately obtain for all  $s \in S_{\infty}$
 \begin{center}
 $ |Q_{s}(a_s x)|_s =  |a_s|_{s}^{2} |Q_{s}(x)|_s \leqslant \varepsilon$ and $ |L_s( a_s x)|_s= |a_s|_{s} |L_{s}(x)|_s  \leqslant \varepsilon$.  
\end{center}
 \smallskip
 Hence given any $\varepsilon > 0$, we get a nonzero element $y = a.x \in  {\mathcal{O}}_{S}^4$ satisfiying the conclusion of Theorem \ref{main thm}, i.e.
\begin{center}
 $ |Q_{s}(y)|_s \leqslant \varepsilon$ and $ |L_s(y)|_s \leqslant \varepsilon$  for all $s \in S$.
\end{center}

%To achieve the proof, one need 
\section{Comments and open problems}~\label{sec end}
\subsubsection*{Irrationality of the pencils forms} The rationality condition in Theorem \ref{main thm}, namely asking irrationality of all the pencils of $Q$ and $L^2$ at all places of $S$ is more restrictive than assuming irrationality of the pencils over $k_S$. Indeed the latter condition leaves the possibility that some pencil could be rational at some place(s) of $S$. In this case, using Ratner's theorem and therefore the classification of intermediate subgroups cannot be achieved by our methods. By analogy with the work of Borel and Prasad in the case of a family of quadratic forms $(Q_s)_{s \in S}$, it may be possible to apply strong approximation and avoiding reduction of dimension. The problem is that this method does not give integral solutions $x\in \mathcal{O}_{S}^n$ of inequalities  $|Q(x)| \leq \varepsilon$ and $|L(x)| \leq \varepsilon$ but only nonzero integral solutions of the pencil forms $|\alpha Q(x) + \beta L^2(x)| \leq \varepsilon$ which may depend on the coefficients $\alpha$ and $\beta$. The most serious issue is to eliminate the dependance on the coefficients, that is, to replace the valid assertion 
$$\forall \varepsilon >0,  \forall P \in \mathbb{P}^{1}(k_S), \exists x \in \mathcal{O}_{S}^{n}-\{0\}, ~~~~~|\widetilde{Q}_P (x)|_S \leq \varepsilon  $$
 by the one we would like
$$\forall \varepsilon >0,  \exists x \in \mathcal{O}_{S}^{n}-\{0\}, \forall P \in \mathbb{P}^{1}(k_S), ~~~~~|\widetilde{Q}_P (x)|_S \leq \varepsilon.  $$
Indeed if one is able to do so,  such $x$ will satisfies those inequalities for both $P_1=[0:1]$ and $P_2=[0:1]$ and by homogeneity it would give the solution of our problem.

\subsubsection*{The problem of null values } The Theorem \ref{main thm} is not conclusive regarding the existence of solutions leading to null values of either $Q$ or $L$. Indeed we are not able to discard this possibility when $s$ is a finite place, the reason is that the use of strong approximation does not provide information enough which guarantees the existence of a solution  $x\in \mathcal{O}_{S}^{n})$ with $0<|Q(x)|_s$ and $<0|L(x)|_s$ for all $s\in S$.

\subsubsection*{Towards density} It should be possible to obtain the density of $F(\mathcal{O}_{S}^{n})$ for a pair $F=(Q,L)$ over $k_S$ under the same assumptions generalising those of Theorem \ref{thm} i.e. without the condition $Q_{| L = 0}$ is nondegenerate added here for our purpose. For this we need a analog of Lemma 6 of \cite{G1} for nonarchimedean completions which has no clear reason to fail in characteristic zero. A significant difference with the classical Oppenheim conjecture is that the stabilizer of such pairs is no more maximal, and the classification of intermediate subgroups is much more involved. Unfortunately we are not able to prove Lemma \ref{archimedean} for non archimedean completions and to avoid the use of strong approximation.

\subsubsection*{An Open problem} We conclude by mentioning a conjecture of Gorodnik (see \cite{G1}, conjecture 15) which concerns the assumption $(2)$ of Theorem \ref{thm} in the real case. It is conjectured that the condition $Q_{| L = 0}$ is isotropic can be replaced by the weaker assumption that the pencil $\alpha Q + \beta L^2$ is isotropic for any real numbers $\alpha,\beta$ such that $(\alpha,\beta) \neq (0,0)$.\\

\begin{conjecture}[Gorodnik]~\label{conj} Let $F  = (Q,L)$ be a pair consisting of one nondegenerate quadratic and one nonzero linear form in dimension $n \geq 4$. Suppose that 
\begin{enumerate}
  \item  For every $\beta \in \bbr$, $Q + \beta L^2$ is indefinite.
  \item For every $(\alpha,\beta) \neq (0,0)$, with $\alpha,\beta \in \bbr$, $\alpha Q + \beta L^2$ is irrational.
 \end{enumerate}
Then $F(\mathcal{P}(\bbz^n))$ is dense in $\bbr^2$.
\end{conjecture}
\noindent The first condition is necessary for density to hold. The main issue is that this condition (contrarily to the condition that $Q_{| L = 0}$ is indefinite) does not allow us to reduce to the four dimensional case. Hence all the strategy of the proof of Theorem \ref{thm} becomes needless regarding the impossibility to classify all the intermediate subgroups in higher dimension.

\end{document}